\documentclass[english,12pt,oneside]{amsproc}
\usepackage[utf8]{inputenc}
\usepackage[T2A]{fontenc}
\usepackage{amsmath,amssymb,amsthm}
\usepackage[a4paper,hmargin=2.5cm,vmargin=2.5cm]{geometry}
\usepackage[english]{babel}
\usepackage{tikz-cd}
\usepackage{enumitem}

\newcounter{stmcounter}[section]
\newcounter{thcounter}

\numberwithin{equation}{section}

\newtheorem{proposition}[stmcounter]{Proposition}
\newtheorem{lemma}[stmcounter]{Lemma}

\newtheorem{theorem}[thcounter]{Theorem}
\newtheorem*{theorem*}{Theorem}

\theoremstyle{definition}
\newtheorem{definition}[stmcounter]{Definition}

\theoremstyle{remark}
\newtheorem{remark}[stmcounter]{Remark}
\newtheorem{example}[stmcounter]{Example}

\newenvironment{pf}{\noindent\textbf{Proof.}}{\qed}
\newcommand{\define}[1]{{\textit{#1}}}

\renewcommand{\leq}{\leqslant}
\renewcommand{\geq}{\geqslant}

\begin{document}

\title{Persistent homological Quillen-McCord theorem}

\author{Vitalii Guzeev}
\address{Faculty of Mathematics, National Research University Higher School of Economics, Russian Federation}
\email{viviag@yandex.ru}

\date{\today}
\thanks{The article was submitted as the graduation thesis of the author at the HSE University.}

\begin{abstract}
The Quillen-McCord theorem (aka Quillen fiber lemma) gives a sufficient condition on a map between classifying spaces of posetal categories to be a homotopy equivalence. Jonathan Ariel Barmak in his paper [J. Comb. Theory Ser. A 118, 8 (November 2011), 2445–2453.] gives an elementary topological proof and proves a homological version of the theorem.

Following his scheme of the proof, we formulate and prove the homological Quillen-McCord theorem, stable with respect to interleaving distances. To formulate the theorem and apply the scheme, we introduce persistence objects as objects in appropriate functor categories, describe a-la barcode decompositions of persistence posets, and prove several results, e.g. order extension principle for objects in Fun(I, Pos) and approximate triviality of left derived functors of approximately trivial objects in Fun(I, R-Mod).

Since the given proof gives explicit Lipschitz constant for the map of persistence classifying spaces, we expect this result to be useful in TDA for reducing complexity of experimental data.
\end{abstract}

\maketitle

\section{Introduction}

Computation of the homotopy type of an arbitrary CW complex is an open problem. It is theoretically possible to compute homology groups in all cases, but in practice, such computations are limited by resources. Complexity generally grows with a dimension of a complex. Hence there is an optimization problem --- given CW complex $B$, construct CW complex $A$ such that $\operatorname{dim}(A) < \operatorname{dim}(B)$ and $H_{\star}(A) \cong H_{\star}(B)$.\\

The Alexandrov-\v{C}ech theorem is proven to be a useful tool in applications [1, 2]. One can associate to a given covering $\mathcal{U}$ a partially ordered set of sets $U \in \mathcal{U}$ and all their intersections, ordered by inclusion. The classifying space of this poset (to be precise, its posetal category) is a barycentric subdivision of geometric realization of the nerve of $\mathcal{U}$, hence these two spaces are homeomorphic. Given this observation, one can formulate the theorem in terms of classifying spaces.\\

Jonathan Barmak proves \cite{Bar11} the homological version of the Quillen-McCord theorem (also known as Quillen fiber lemma or Quillen's theorem A for posets). It can be stated as follows.
\begin{theorem*}
  Assume $X, Y$ are finite posets, $f : X \to Y$ is an order-preserving map, $R$ is a PID.\\
  If $H_i(\mathcal{B}(f^{-1}(Y_{\leqslant y})),R) = 0$ for any $i$ and $y \in Y$, then $\mathcal{B}f$ induces isomorphisms of all homology groups with coefficients in $R$ on $BX$ and $BY$.
\end{theorem*}

Provided an algorithm for the construction of $X$ and a map by $Y$, this theorem may provide a partial solution to the stated optimization problem, in particular, for nerves of coverings.\\

Coverings we operate may come from a series of observations, for instance, from a long experiment. While considering the series as a whole, we are interested in the persistence of nerves. To operate experimental data some stable version of the theorem is required.\\

In this paper, we formulate and prove the persistent homological version of the Quillen-McCord theorem.
\begin{theorem*}
    Assume $X, Y$ are persistence posets of finite type indexed by a very good monoid $I$, $f : X \to Y$ is an order-preserving map. Let $m$ be the number of elements of $Y$ and $R$ be a PID. Assume that the persistence topological space $\mathcal{B}(f^{-1}(Y_{\leqslant y}))$ is $\varepsilon$-acyclic over $R$ for any $y=(\ldots,y_i,\ldots) \in Y$. Then $BX$ and $BY$ are $4m\varepsilon$-interleaved over $R$.\\
\end{theorem*}

The terms used in the statement are explained in sections 2 and 3.\\

The paper is structured as follows.
\begin{enumerate}
  \setcounter{enumi}{1}
  \item In preliminaries we give an outline of well-known notions used throughout the paper.
  \begin{enumerate}[label*=\arabic*.]
    \item First preliminary subsection is devoted to the notion of the interleaving distance (following {\cite{GS16}}) which fits as a required measure of similarity between series. Usability of this notion is guaranteed by persistence theorem {\cite{Zomorodian05}}. We try to keep some level of generality, in particular, we work with a generalized version of the theorem proven in {\cite{Corbet18}} and formulate it. This generality is significant for applications and allows us to operate both discrete and continuous models of time in an experiment.
    \item Second subsection gives a classical formulation of the Quillen-McCord theorem and necessary definitions. To bind the theorem to an abstract context in which the general Quillen A theorem {\cite{Quillen72}} is used, we derive it from the general theorem. Additionally, we outline Barmak's ({\cite{Bar11}}) proofs of classical and homological versions of the theorem.
  \end{enumerate}
  \item In the apparatus section we systematize and sometimes introduce various results forming the toolchain to prove the target theorem. Its necessary set of definitions is close to well-established in a field (for instance, {\cite{Bubenik15}}), but the whole toolchain is independent.
  \begin{enumerate}[label*=\arabic*.]
    \item At first, we reformulate several preliminary definitions in terms of functor categories and transfer some notions, e.g. classifying space of a category and a covering, to appropriate functor categories.
    \item In the second section we formulate and prove the linear extension principle for functors to the category of posets. We give both finite and general versions with an AC-dependent proof of the latter.
    \item At the end of the section we give some technical stability results concerning objects of appropriate functor categories.
  \end{enumerate}
  \item Finally, we formulate and prove the main result --- stable (w.r.t. interleaving distances) homological Quillen-McCord theorem. We adapt a proof given in ({\cite{Bar11}}) using the developed apparatus.
\end{enumerate}

\section{Preliminaries}

\subsection{Persistence modules and interleaving distance}

We start with the definition of a simplicial set.

\begin{definition}
  The simplex category $\Delta$ is a category of non-empty totally-ordered sets of finite length with order-preserving functions as morphisms.
\end{definition}

\begin{definition}
  \define{Simplicial set} is a contravariant functor from $\Delta$ to $Set$.
\end{definition}

\begin{proposition}
  Simplicial sets form the category $sSet$.
\end{proposition}

The real initial definition, representing a series of simplicial sets, is the following.

\begin{definition}
  \define{Persistence simplicial set} is a family of simplicial sets $C_{0} \xrightarrow{f_0} C_{1} \xrightarrow{f_1} C_{2} \xrightarrow{f_2} \ldots$ where $f_i$ are natural transformations. We call the maps $f=(\ldots,f_i,\ldots)$ the \define{structure maps} of a persistence simplicial set.
\end{definition}

\begin{definition}
  Let $R$ be a ring. \define{Persistence module over $R$} is a family of $R$-modules $M_i$ with homomorphisms $\phi_i : M^i \to M^{i+1}$ as the structure maps. We denote the composition of structure maps between $M_i$ and $M_j$ by $\phi_{ij}$.
\end{definition}

An example of a persistence module is given by homology modules of persistence simplicial set $C_{\star}$ (\define{persistent homology}). We set $H_i^j(C_{\star}) := H_i(C_{j})$, the maps $\phi_j$ are induced by $f_i$.\\

\begin{definition}
  Persistence simplicial set (module) is of \define{finite type} over $R$ if all its simplicial sets (modules) have a finite number of non-empty images (are finitely generated as $R$-modules) and all $f_i$ ($\phi_i$) are isomorphisms for $i > m$ for some $m$.
\end{definition}

\begin{definition}
  A persistence module is of \define{finitely presented type} over $R$ if all its modules are finitely presented as $R$-modules and all $f_i$ ($\phi_i$) are isomorphisms for $i > m$ for some $m$.
\end{definition}

We shall generally use the terms of modules in this section.\\

One can study maps of persistence modules. 

\begin{definition}
  Let $M$ and $N$ be persistence modules. The family $f$ of maps $f_i : M_i \to N_i$ is called a morphism from $M$ to $N$ if all $f_i$s commute with structure maps.
\end{definition}

\begin{definition}
  Let $M$ and $N$ be persistence modules, $f$ is a collection of maps $f_i : M_i \to N_{i+\varepsilon}$. Then if all $f_i$ commute with structure maps, $f$ is called an \define{$\varepsilon$-morphism}.
\end{definition}

There is a general notion of interleaving distance between persistence modules.

\begin{definition}
  We denote by $\operatorname{Id}_{\varepsilon} : M \to M$ the shift of persistence module, defined by compositions of structure maps $M_i \to M_{i+\varepsilon}$ for all $i$.
\end{definition}

\begin{definition}
  Persistence modules $M$ and $N$ are called \define{$\varepsilon$-interleaved} ($M \stackrel{\varepsilon}{\sim} N$) if there exists a pair of $\varepsilon$-morphisms $(\phi : M \to N,\;\psi : N \to M)$ called an \define{$\varepsilon$-interleaving} such that $\phi \circ \psi = \operatorname{Id}_{2\varepsilon} : N \to N$ and $\psi \circ \phi = \operatorname{Id}_{2\varepsilon} : M \to M$.\\
\end{definition}

\begin{remark}
  It follows that $M \stackrel{\varepsilon}{\sim} N$ implies $M \stackrel{\alpha}{\sim} N$ for any $\alpha > \varepsilon$ since for an $\varepsilon$-interleaving $(\phi, \psi)$ we have an $\alpha$-interleaving $(\operatorname{Id}_{\alpha - \varepsilon} \circ \phi, \operatorname{Id}_{\alpha - \varepsilon} \circ \psi)$.
\end{remark}

\begin{definition}
  The $\varepsilon$-interleaving induces an extended pseudometric on a set of persistence modules. This pseudometric is defined as $d(X,Y) = min\{\varepsilon \in I\;|\;X \stackrel{\varepsilon}{\sim} Y\}$. This pseudometric is called interleaving distance. {\cite[Definition 2.12]{GS16}}\\
\end{definition}

There is a well-known theorem.

\begin{theorem} {\cite[Theorem 3.1]{Zomorodian05}}\\
  The category of persistence modules of finite type over Noetherian ring with unity $R$ is equivalent to the category of graded finitely generated $R[t]$-modules.
\end{theorem}

It is proven in {\cite{Corbet18}}. The authors provide a generalization that is more suitable for our needs.\\

\begin{definition}
  Let $(G,\star)$ be a commutative monoid and $g_1, g_2 \in G$.
  We say that $g_1 \preceq g_2$ if $\exists h \in G:\; h \star g_1 = g_2$, $h$ is not a neutral element.
\end{definition}

\begin{definition}
  Let $(G,\star)$ be a monoid. Consider $H \subset G$. An element $m \in G$ is called a \define{common multiple} of $H$ if $h \preceq m$ for any $h \in H$. A common multiple $m$ of $H$ is called \define{partially least} if there is no common multiple $m_1$ of $H$ such that $m_1 \preceq m$.
\end{definition}

\begin{definition} {following \cite[Definition 11]{Corbet18}}
  Monoid $(G,\star)$ is called \define{good} if the following hold:
  \begin{itemize}
    \item $(G, \star)$ is commutative;
    \item $g_1 \star g_2 = g_1 \star g_3$ implies $g_2 = g_3$ (cancellation);
    \item $g_1 \preceq g_2$ and $g_2 \preceq g_1$ imply $g_1 = g_2$ (anti-symmetricity);
    \item For any finite $H \subseteq G$ there exists at most finitely many partially least common multiples (property of being weak plcm).
  \end{itemize}
\end{definition}

\begin{definition} {\cite[Definition 12]{Corbet18}}
  Let $R$ be a ring and $G$ be a good monoid. \define{(Generalized) persistence module} is a family of $R$-modules $M^i$ for $i \in G$ with homomorphisms $\phi_{ij} : M^i \to M^j$ for $i \preceq j$, satisfying identity relation $\phi_{ii} = \operatorname{Id}$ and composition relation $\phi_{ij} \circ \phi_{jk} = \phi_{ik}$ for any $i \preceq j \preceq k \in G$ as the structure maps.
\end{definition}

\begin{theorem} {\cite[Theorem 21]{Corbet18}}\\
  Let $R$ be a ring with unity and $G$ be a good monoid. Then the category of finitely presented graded $R[G]$-modules is isomorphic to the category of $G$-indexed persistence modules over $R$ of finitely presented type.
\end{theorem}

\begin{remark}
  Cancellative commutative monoid with $\preceq$ being a total order (totally ordered) is good.
\end{remark}

\begin{definition}
  We call cancellative commutative totally ordered monoids \define{very good}.
\end{definition}

\begin{example}
  Consider the monoid of non-negative real numbers $\mathbb{R}_{\geq 0}$ with addition as a monoidal operation. The addition of real numbers is cancellative and $\preceq$ is a total order. This monoid is very good.\\

  For contrast, we can consider $\mathbb{R}$. Since $g_1 = h + g_2$ implies $g_2 = (-h) + g_1$ $\preceq$ is trivial and $R$ is not a good monoid in the sense of the given definition.
\end{example} ~ \par

We can examine how notions related to interleaving distance look in a category of graded modules under additional constraints on the indexing set.

\begin{definition} {\cite[Definition 2.7]{GS16}}\\
  Let $M$ and $N$ be graded $R[G]$-modules, $f : M \to N$ be a homomorphism of modules. Then $f$ is called \define{$\varepsilon$-morphism} if $f(M^j) \subset N^{j+\varepsilon}$.
\end{definition}

Let $R$ be a ring with unity and $G$ be a very good monoid. Then for any graded $R[G]$-module $M$ and any $\varepsilon \in \mathbb{G}$ there exists a $\varepsilon$-morphism $\operatorname{Id}_{\varepsilon} : M \to M$. For instance, if the monoid $G$ is the monoid of non-negative integers, the $\operatorname{Id}_{\varepsilon}$ is a multiplication by $t^{\varepsilon}$, where $t$ is a generator of a polynomial algebra {\cite[Example 2.8]{GS16}}.\\

This example can be generalized as follows.

\begin{proposition} \cite[Following the Equation 4]{Corbet18}
  \label{mult}
  Let $M$ be an $R[G]$-module. The $\varepsilon$-morphism $\operatorname{Id}_{\varepsilon} : M \to M$ is a multiplication by a fixed element $m \in R[G]^{\varepsilon}$, where grading on $R[G]$ is given by $G$.
\end{proposition}

\begin{proposition} {\cite[Proposition 2.13]{GS16}}
  \label{epstriv}
  Condition $M \stackrel{\varepsilon}{\sim} 0$ is equivalent to condition $m^2M = 0$.
\end{proposition}

\begin{lemma}
  \label{ops}
  Let $0 \to M \to L \to N \to 0$ be a short exact sequence of graded modules. Then the following properties hold.
  \begin{itemize}
    \item If $M \stackrel{\varepsilon_1}{\sim} 0$ and $N \stackrel{\varepsilon_2}{\sim} 0$ then $L \stackrel{\varepsilon_1 + \varepsilon_2}{\sim} 0$. {\cite[Proposition 4.6]{GS16}}
    \item If $L \stackrel{\varepsilon}{\sim} 0$ then $M \stackrel{\varepsilon}{\sim} 0$ and $N \stackrel{\varepsilon}{\sim} 0$.
    \item If $M \stackrel{\varepsilon}{\sim} 0$ then $L \stackrel{2\varepsilon}{\sim} N$. {\cite[Proposition 4.1]{GS16}}
    \item If $N \stackrel{\varepsilon}{\sim} 0$ then $M \stackrel{2\varepsilon}{\sim} L$. {\cite[Proposition 4.1]{GS16}}
  \end{itemize}
\end{lemma}

The second statement of the lemma requires proof. We give it for non-negative integers first for compatibility with the referenced result. 

\begin{pf}
  Denote non-trivial maps in s.e.s as $i$ and $q$.\\

  Then $i(t^{2\varepsilon}a) = t^{2\varepsilon}i(a) = 0$ for any $a \in M$. The homomorphism $i$ is injective. Hence by Proposition \ref{epstriv} we have $M \stackrel{\varepsilon}{\sim} 0$.\\

  On the other side, we have $0 = q(t^{2\varepsilon}a) = t^{2\varepsilon}q(a)$ where $q$ is surjective. Hence $N \stackrel{\varepsilon}{\sim} 0$.
\end{pf}\\

Proposition \ref{mult} allows to generalise the lemma by changing $t^{\varepsilon}$ to $m$.\\

By equivalence of categories, the lemma holds for persistence modules over rings with unity indexed by very good monoids.

\subsection{The Quillen-McCord theorem}

\begin{definition}
  Let $X$ be a set and $\mathcal{P}(X)$ be its powerset. A set $S \subset \mathcal{P}(X)$ is called a simplicial complex if $W \in S$ for any $V \in S$ and $W \subset V$. An element of a simplicial complex of cardinality $n$ is called $n$-simplex.
\end{definition}

\begin{proposition}
  Simplicial complexes form the category $sCpx$.
\end{proposition}

\begin{definition}
  The \define{join} $A \star B$ of simplicial complexes $A$ and $B$ is the simplicial complex with simplices --- all possible unions of simplices $a \in A$ and $b \in B$.
\end{definition}

Let $A$ be a simplicial complex.

\begin{definition}
  The \define{star} $\operatorname{st}(x)$ of simplex $x \in A$ is the minimal by inclusion simplicial complex containing all simplices $a \in A$ such that there exists an inclusion $x \hookrightarrow a$.
\end{definition}

\begin{definition}
  \define{Link} $\operatorname{lk}(x)$ of simplex $x \in A$ is defined as follows: $\operatorname{lk}(x) = \{v \in \operatorname{st}(x)|\; x \cap v = \varnothing\}$.
\end{definition}

\begin{proposition}
  For $0$-simplex $x$ it holds that $\operatorname{st}(x) = \operatorname{lk}(x) \star x$.
\end{proposition}

\begin{definition}
  Functor $\left[\bullet\right] : sCpx \to Top$ which maps $n$-simplices to standard geometric $n$-simplices is called the standard \define{geometric realization} of a simplicial complex.
\end{definition}

\begin{definition}
  Functor $\left|\bullet\right| : sSet \to Top$ which maps $n$-simplices to $n$-cells and face maps to attachment maps is called the \define{geometric realization} of a simplicial set.
\end{definition}

\begin{definition}
  \define{Join} of topological spaces $A$ and $B$ is defined as follows: $A \star B\ := A \sqcup_{p_0} (A \times B \times [0,1]) \sqcup_{p_1} B$, where $p$ are projections of the cylinder $A \times B \times [0,1]$ onto faces.
\end{definition}

The next proposition gives motivation for the definition of a join of simplicial complexes.

\begin{proposition}
  $\left|A \star B\right| = \left|A\right| \star \left|B\right|$. Hence if $x$ is a $0$-simplex, then $\left|\operatorname{st}(x)\right|$ is a cone over $\left|\operatorname{lk}(x)\right|$.
\end{proposition}

\begin{definition}
  There is a standard set of definitions, accompanying the definition of a simplicial set. Let $S$ be a simplicial set. Then
  \begin{enumerate}
    \item images $S([n])$ of sets of cardinality $n$ are called $n$-simplices;
    \item images of injective maps $S([n] \to [n+1])$ are called face maps;
    \item images of surjective maps $S([n+1] \to [n]$ are called degeneracy maps.
  \end{enumerate}
  
  Details of the explicit construction can be found in \cite{Friedman}
\end{definition}

\begin{definition}
  Let $\mathcal{C}$ be a small category. Then we can functorially (w.r.t. to category $Cat$) assign to it a simplicial set $\mathcal{N}(\mathcal{C})$ called \define{nerve of a category} to it.\\

  Construction goes as follows:
  \begin{itemize}
    \item we assign to each object of $\mathcal{C}$ a $0$-simplex and to each morphism in $\mathcal{C}$ a $1$-simplex with order following corresponding arrow in a category;
    \item then we take the set of all morphisms as an alphabet and write all the words in it such that we can move from the first letter to the last following arrows in a category. We assign an $l$-simplex to a word of length $l$. Each commutative triangle $f, g, h = f \circ g$ in $\mathcal{C}$ gives rise to a morphism between these words --- it replaces $f \circ g$ with $h$ and serves as a face map. Replacements $f \to f \circ Id$ represent degeneracy maps.
  \end{itemize}

  This collection of data forms a simplicial set, details can be found in \cite{Kerodon002M}.
\end{definition}

\begin{definition}
  Geometric realization $B\mathcal{C}$ of $\mathcal{N}(\mathcal{C})$ is called \define{classifying space} of $\mathcal{C}$.\\

  We denote the composition of nerve and geometric realization as $\mathcal{B}$. It is a composition of functors, hence a functor. By definition $\mathcal{B}(\mathcal{C}) = B\mathcal{C}$ and we prefer notation $\mathcal{B}(f)$ to $Bf$.
\end{definition}

\begin{definition}
  Let $F: \mathcal{C} \to \mathcal{D}$ be a functor and $d$ --- object in $\mathcal{D}$. Then \define{comma category} $d \downarrow F$ is a category with objects --- pairs $(s,i_s)$ of objects in $\mathcal{C}$ and morphisms $i_s : d \to F(s)$ and morphisms --- morphisms $g$ in $\mathcal{C}$ such that triangle $i_s, i_{g(s)}, F(g)$ is commutative.
\end{definition}

\begin{theorem} {\cite[Theorem A]{Quillen72}}\\
  If $F: C \to D$ is a functor such that the classifying space $B(d \downarrow F)$ of the comma category $d \downarrow F$ is contractible for any object $d \in D$, then $F$ induces a homotopy equivalence $BC \to BD$.
\end{theorem}

The nerve construction on a posetal category yields a simplicial set with the following good properties: the image of a set of reorderings on any set in $\Delta$ contains unique simplex and subsets of its simplices are also simplices. These simplicial sets form a subcategory $ssCpx$ of $sSet$.\\

There is a functor $F$ from this subcategory to the category of simplicial complexes, that turns a simplicial set to a set of its images.\\

\begin{proposition}
  Let $X$ be an object of $ssCpx$. Then $|X|$ is homeomorphic to the geometric realization of simplicial complex $[F(X)]$.
\end{proposition}

This proposition allows us to operate simplicial complexes instead of simplicial sets and avoid general definitions of join, star, and link for simplicial sets.\\

Application of Quillen's A theorem to posets yields the following theorem (we identify poset with its posetal category).

\begin{theorem} \textbf{the Quillen-McCord theorem}\\
  Assume $X, Y$ are finite posets, $f : X \to Y$ is an order-preserving map.\\
  If $\forall y \in Y\;\mathcal{B}(f^{-1}(Y_{\leqslant y}))$ is contractible, then $\mathcal{B}f$ is a homotopy equivalence between $BX$ and $BY$.\\
\end{theorem}

\begin{theorem} \textbf{Homological Quillen-McCord theorem} \cite[Corollary 5.5]{Bar11}\\
  Assume $X, Y$ are finite posets, $f : X \to Y$ is an order-preserving map, $R$ is a PID.\\
  If $\forall y \in Y\;H_i(\mathcal{B}(f^{-1}(Y_{\leqslant y})),R) = 0$ for any $i$, $\mathcal{B}f$ induces isomorphisms of all homology groups of $BX$ and $BY$ with coefficients in $R$.\\
\end{theorem}

Proofs of both theorems are used in our arguments and we recall them in brief.

\subsubsection{Barmak's proof of the Quillen-McCord theorem}

\begin{proposition} {Variation of \cite[Lemma 2.2]{Bar11}}
  \label{prop:comparison}
  Let $f,g : X \to Y$ be order-preserving maps between finite posets such that $\forall x\;f(x) \leq g(x)$. Then $\mathcal{B}(f)$ is homotopy-equivalent to $\mathcal{B}(g)$.
\end{proposition}

\begin{proposition}
  Note that $\operatorname{lk}(F(\mathcal{N}(x))) = F(\mathcal{N}(X_{>x})) \star F(\mathcal{N}(X_{<x}))$. Therefore $\left|\operatorname{lk}(\mathcal{N}(x))\right| = \mathcal{B}(X_{>x}) \star \mathcal{B}(X_{<x})$.
\end{proposition}

The crucial observation is the existence of the following covering for any $x \in X$.

\begin{equation}\label{covering}
  \mathcal{B}(X) = \mathcal{B}(X \setminus \{x\}) \cup \left|\operatorname{st}(\mathcal{N}(x))\right|.\\
\end{equation}

\begin{lemma}
  \label{lem:homotopy}
  Let $X$ be a finite poset and for $x \in X$ either $\mathcal{B}(X_{>x})$ or $\mathcal{B}(X_{<x})$ is contractible. Then embedding $\mathcal{B}(X \setminus \{x\}) \hookrightarrow \mathcal{B}(X)$ is a homotopy equivalence.
\end{lemma}

For this paper, this lemma is the most important part of the proof of the theorem, hence we recall its proof.\\

\begin{pf}
  By hypothesis, the space $\left|\operatorname{lk}(\mathcal{N}(x))\right| = \left|\operatorname{st}(\mathcal{N}(x))\right| \cap \mathcal{B}(X \setminus \{x\})$ is contractible. Hence its embedding to its cone $\left|\operatorname{st}(\mathcal{N}(x))\right|$ is a homotopy equivalence by Whitehead theorem. Being a subcomplex, it is a strong deformation retract. Then, $\mathcal{B}(X \setminus \{x\})$ is a strong deformation retract of $\mathcal{B}X = \left|\operatorname{st}(\mathcal{N}(x))\right| \cup \mathcal{B}(X \setminus \{x\})$
\end{pf}

\begin{definition} {Variation of \cite[Proposition 2.1]{Bar11}}
  Let $f : X \to Y$ be an order-preserving map between posets. Denote orders ($\leq$) on $X$ and $Y$ as $R_X$ and $R_Y$. Then we define poset $M(f) = X \coprod_f Y$ with $R = R_X \cup R_Y \cup R_{f}$ where $(x,y) \in R_f$ if and only if $(f(x),y) \in R_Y$.\\

  We call this poset a \define{mapping cylinder} of $f$. There are also defined canonical inсlusions $i_X : X \to M(f)$ and $i_Y : Y \to M(f)$.
\end{definition}

\begin{pf} \textbf{the Quillen-McCord theorem}\\
  Let $X, Y$ be finite posets with an order-preserving map $f : X \to Y$.\\

  Every poset has a linear extension. Let $x_1, x_2, \ldots, x_n$ be the enumeration of $X$ in an arbitrary linear extension and $Y^r = \{x_1,\ldots,x_r\} \cup Y \subset M(f)$ for any $r$.\\

  Consider $Y^r_{>x_r} = Y_{\geq f(x_r)}$. The space $\mathcal{B}(Y_{\geq f(x_r)})$ is a cone over $\mathcal{B}(f(x_r))$. It is contractible, therefore $\mathcal{B}(Y^{r-1}) \hookrightarrow \mathcal{B}(Y^{r})$ is a homotopy equivalence by Lemma \ref{lem:homotopy}. By iteration, $\mathcal{B}(j) : \mathcal{B}(Y^{0}) = \mathcal{B}(Y) \hookrightarrow \mathcal{B}(M(f)) = \mathcal{B}(Y^n)$ is homotopy equivalence between $BY$ and $M(f)$.\\

  Then consider a linear extension of $Y$ with enumeration $y_1,\ldots,y_m$ and $X^r = X \cup \{y_{r+1},\ldots,y_m\} \subset M(f)$. We have $X^{r-1}_{< y_r} = f^{-1}(Y_{\leqslant y_r})$. The classifying space of the latter is contractible by the assumption of the theorem. Hence $\mathcal{B}(X^{r}) \hookrightarrow \mathcal{B}(X^{r-1})$ is a homotopy equivalence and by transitivity $\mathcal{B}(i_X)$ is a homotopy equivalence between $X$ and $M(f)$.\\

  Note that $i(x) \leqslant (i_Y \circ f)(x)$. By Proposition \ref{prop:comparison} the space $\mathcal{B}(i_X)(X)$ is homotopy-equivalent to the space $\mathcal{B}(i_Y \circ f)(X) = (\mathcal{B}(i_Y) \circ \mathcal{B}(f))(X)$. Hence $\mathcal{B}(f)$ is the homotopy equivalence between $BX$ and $BY$.
\end{pf}

\subsubsection{Barmak's proof of the homological Quillen-McCord theorem}

\begin{proposition} {\cite[Lemma 2.1]{Milnor56}}
  \label{prop:kunneth}
  Reduced homology modules with coefficients in a principal ideal domain $R$ satisfy the relation
  $H_{r+1}(A \star B, R) \simeq \bigoplus_{i+j=r}(H_i(A,R) \otimes_R H_j(B,R)) \oplus \bigoplus_{i+j=r-1} \operatorname{Tor}_1^R(H_i(A,R),H_j(B,R))$.
\end{proposition}

\begin{lemma}
  \label{lem:trivialties}
  Let $X$ be a finite poset and for $x \in X$ either $H_i(\mathcal{B}(X_{< x}))$ or $H_i(\mathcal{B}(X_{> x}))$ with coefficients in a PID are equal to the homology of a point. Then embedding $\mathcal{B}(X \setminus \{x\}) \hookrightarrow \mathcal{B}(X)$ induces isomorphisms of all homology groups.
\end{lemma}

\begin{pf} ~ \par
  By Proposition \ref{prop:kunneth}, $H_i(\left|\operatorname{lk}(\mathcal{N}(x))\right|) = H_i(\mathcal{B}(X_{>x}) \star \mathcal{B}(X_{<x}))$ are trivial for all indices $i$.
  Application of Mayer-Vietoris long exact sequence to covering introduced in equation \eqref{covering} yields the lemma.
\end{pf}\\

Proof of the theorem is similar to the proof at the end of the previous subsection. We write it here to highlight differences. Changed parts are written in italic.\\

\begin{pf} \textbf{the homological Quillen-McCord theorem}\\
Let $X, Y$ be finite posets with an order-preserving map $f : X \to Y$.\\

Every poset has a linear extension. Let $x_1, x_2, \ldots, x_n$ be the enumeration of $X$ in a fixed linear extension and $Y^r = \{x_1,\ldots,x_r\} \cup Y \subset M(f)$ for any $r$.\\

Consider $Y^r_{>x_r} = Y_{\geq f(x_r)}$. The space $\mathcal{B}(Y_{\geq f(x_r)})$ is a cone over $\mathcal{B}(f(x_r))$. It is contractible, therefore $\mathcal{B}(Y^{r-1}) \hookrightarrow \mathcal{B}(Y^{r})$ is a homotopy equivalence by Lemma \ref{lem:homotopy}. By iteration, the map $\mathcal{B}(j) : \mathcal{B}(Y^{0}) = \mathcal{B}(Y) \hookrightarrow \mathcal{B}(M(f)) = \mathcal{B}(Y^n)$ is a homotopy equivalence between $BY$ and $M(f)$.\\

Now consider a linear extension of $Y$ with enumeration $y_1,\ldots,y_m$ and $X^r = X \cup \{y_{r+1},\ldots,y_n\} \subset M(f)$. We have $X^{r-1}_{< y_r} = f^{-1}(Y_{\leqslant y_r})$. \textit{The classifying space of the latter is acyclic over $R$ by the assumption of the theorem. Hence $\mathcal{B}(X^{r}) \hookrightarrow \mathcal{B}(X^{r-1})$ induces isomorphisms of all homology groups and by the functoriality of homology $\mathcal{B}(i)$ induces isomorphisms of all homology groups between $X$ and $M(f)$.}\\

Note that $i(x) \leqslant (j \circ f)(x)$. By Proposition \ref{prop:comparison} $\mathcal{B}(i)$ is homotopic to $\mathcal{B}(j \circ f) = \mathcal{B}(j) \circ \mathcal{B}(f)$. \textit{Homotopic maps induce the same maps of homology modules, $j$ is a homotopy equivalence and induces isomorphisms. Hence $\mathcal{B}(f)$ induces isomorphisms between $H_i(BX,R)$ and $H_i(BY,R)$.}
\end{pf}\\

We see two updates. The first one is essential, it requires Lemma \ref{lem:trivialties} and operates some equivalence propagating in a chain of length equal to the cardinality of $Y$. The second follows automatically from the functoriality of all used constructions.

\section{Apparatus for main results}

\subsection{Persistence objects and related constructions}

We have two types of persistence objects with similar definitions. There is a general notion of a persistence object such that these definitions fall into special cases.\\

\begin{definition}
  Consider $I$ --- the poset category of a fixed linearly ordered set. There is a sequence category $Fun(I, \mathcal{C})$ of functors from $I$ to some category $\mathcal{C}$. We call objects of this category \define{persistence objects over $\mathcal{C}$}
\end{definition}

\begin{example} ~ \par
  \begin{itemize}
    \item Persistence complex is a persistence object over the category of chain complexes;
    \item Persistence $R$-module is a persistence object over the category of $R$-modules; \cite[Definition 2.15]{GS16}
    \item Persistence simplicial set is a persistence object over the category $sSet$;
    \item Persistence poset is a persistence object over $Pos$;
    \item Persistence topological space is a persistence object over $Top$.
  \end{itemize}
\end{example}

\begin{definition}
  We call images of morphisms in $I$ the \define{structure maps} of a persistence object over $\mathcal{C}$.
\end{definition}

We use the notation $(X,\phi)$ for ``Persistence object X with the family of structure maps $\phi$ over fixed indexing category $I$''. We use notation $\phi_{ij}$ for a structure map between $X_i$ and $X_{j}$.\\

Consider $\mathcal{F}$ --- a functor from $\mathcal{C}$ to $\mathcal{D}$. It naturally extends to a functor between $Fun(I,\mathcal{C})$ and $Fun(I,\mathcal{D})$. Let $P$ be a persistence poset. Apparently, $\mathcal{B}(P)$ is a persistence topological space.\\

\begin{definition}
  A persistence topological space $X$ is called \define{$\varepsilon$-acyclic} over $R$ if for all indices $H_i(X,R) \stackrel{\varepsilon}{\sim} H_i(pt,R)$.
\end{definition}

\begin{definition}
  A finite sequence of finite posets is called a \define{persistence poset of finite type}.
\end{definition}

\begin{remark}
  The definition of a persistence poset has a foundational issue --- since there are no maps to an empty set, structure maps are in general not defined. In particular, this issue is unavoidable in the case of posets of finite type. We prefer to resolve this issue by adjoining a contradiction symbol $\bot$ to a class of functions. If there is no candidate for a function $\phi_{ij}$ in $Pos$, we set $\phi_{ij} = \bot$ and ignore these functions in reasonings similar to ``something is defined by all images of $x$ under the structure maps''.
\end{remark}

\begin{proposition}
  Let $X$ be a persistence poset of finite type. Then $\mathcal{B}(X)$ has homology modules of finitely presented type.
\end{proposition}

\begin{pf}
  The nerve of an empty poset is an empty simplicial set. Hence $\mathcal{N}(X)$ is a finite sequence of simplicial sets. Since a finite poset has only a finite number of chains, each component is finite, hence $\mathcal{N}(X)$ has finite type. The homology of $\mathcal{B}(X)$ can be computed as the homology of this simplicial set. Hence homology modules of $\mathcal{B}(X)$ are of finitely presented type as quotients of finite rank free modules by finite rank free modules.
\end{pf}\\

\begin{definition}
  Let $(X, \phi)$ and $(Y, \psi)$ be persistence posets and $f : X \to Y$ be a morphism. Consider the series of mapping cylinders of posets $M(f_i)$. We can define the structure maps $\chi$ on this series as $\phi$ on $X$ and $\psi$ on $Y$. They preserve orders $R_{f_i}$ by definition of structure maps. The persistence poset $(M(f), \chi_i)$ is a \define{mapping cylinder of a map of persistence posets}. We also have canonical inclusions $i_X$ and $i_Y$.
\end{definition}

We can also define a subobject in $Fun(I, Pos)$.

\begin{definition} ~ \par
  Consider persistence poset $(X,\phi)$ and the collection of subsets $Y_i \subset X_i$ for $i \in I$.
  If for any $i \in I$, $y \in Y_i$ and for any $j>i$ we have $\phi_{ij}(y) \in Y_j$ then $(Y,\phi)$ is called a \define{persistence subposet} of $(X,\phi)$.
\end{definition}

We can also define an \define{element} of a persistence poset.
\begin{definition}
  A persistence subposet $x$ of $(X,\phi)$ with at most one element in each component is called an element of $X$.
\end{definition}

\begin{remark}
  Consider the persistent poset $(X, \phi)$ of finite type. For each element $x \in X$, there exists a minimal index $i$ such that $x_i \neq \emptyset$. Then $x$ is completely defined by its images under structure maps.
\end{remark}

One can consider component-wise order on $X$ and work with sets like the following.
\begin{example}
  Let $x$ be an element of the persistence subposet $(X, \phi)$. Then $X_{<x}$ is the poset of elements component-wise less than $x$.\\

  $X_{<x}$ contains component-wise only elements comparable with $x$. Since structure maps are order-preserving, $X_{<x}$ is a persistence subposet.
\end{example} ~ \par

Finally, we give a definition of \define{persistence covering}.\\

\begin{definition}
  Assume a poset $X$ splits into a union of subposets $X_j$. Then every subposet $X_j$ has its own classifying space $BX_j$. If these spaces (or minimal open sets containing them in $BX$) cover the whole $BX$, they are called a \define{persistence covering}.
\end{definition}

This definition gives an example of how structures in the category of persistence posets can be transferred to other persistence categories. It is possible to reformulate the definition as internal to the category of persistence topological spaces but we prefer to keep a more constructive way.

\subsection{Order extension principle for persistence posets}

In his proof of the Quillen-McCord theorem Barmak relies on the order extension principle. To be able to transfer Barmak's proof to the persistent case we have to stress a similar statement for persistence posets.

\begin{definition}
  An \define{extension of persistence poset} $X$ is a series of partially-ordered extensions of $X_i$ such that the structure maps of $X$ are well-defined on these extensions. If extensions of all components are linear, we call this series a linear extension.
\end{definition}

\begin{proposition}
  \textbf{Transfer of order}. Let $f$ be a morphism between posets $X$ and $Y$ and $\overline{Y}$ be a linear extension of $Y$. Then $f$ induces partially ordered extension $\hat{X}$ of $X$ such that $f$ is well-defined as map $\hat{X} \to \overline{Y}$.
\end{proposition}

\begin{pf}
  Consider two incomparable points $a, b \in X$ and map $f : X \to \overline{Y}$ which is well-defined. One of the following holds.\\
  \begin{itemize}
    \item f(b) < f(a)
    \item f(a) < f(b)
    \item f(a) = f(b)
  \end{itemize}
  If strict inequality holds, we can impose a single relation on $a$ and $b$ --- we inherit relation from images.\\
  If equality holds we do not add any new relation.
\end{pf}

\begin{proposition}
  \textbf{Left propagation of linear extension.} Assume indexing set $I$ of persistence poset $(X,\phi)$ is converse well-founded and there exists maximal index $i$ such that $X_j = \emptyset$ for any $j > i$. Then $(X,\phi)$ has a linear extension.
\end{proposition}

\begin{pf}
  We can extend the order on the component $X_i$ to linear. Given this order, we can transfer it to the left via all structure maps. We obtain an extension of $(X,\phi)$ because all preimages of incomparable elements were incomparable and we have equipped them with compatible orders. Now let's assume we obtained linear orders in components $X_j$ for $j > j_0$ by this construction. Then $X_{j_0}$ can be linearly extended. The proposition follows by transfinite induction and by simple induction if $X$ is finite as a sequence.
\end{pf}\\

This statement can also be seen as a corollary of a more general proposition.
One can consider set $E(X,\phi)$ of extensions of persistence poset $(X,\phi)$ with partial order defined as follows: let $Y, Z$ be extensions of $X$, then $Y \geq Z$ if and only if $Y$ is an extension of $Z$. Since underlying sets of these extensions are always the same we can identify elements of $E(X,\phi)$ with tuples of order relations on components of $X$.

\begin{proposition}
  \label{Zorn_cond}
  Every linearly ordered subset of $E(X,\phi)$ has an upper bound in $E(X,\phi)$.
\end{proposition}

\begin{pf}
  Let $\{R^s|\;s \in S\}$ be a linearly ordered subset of $E(X,\phi)$ indexed by set $S$. Consider $R = \bigcup R^s$ where the union is taken component-wise. Assume for some elements $a, b \in X_i$ for some $i$ some $\phi_{ij}$ cannot be defined on them as an order-preserving map. Then there exists $s \in S$ such that $a$ and $b$ are comparable in extension $R^s$. But $R^s$ is an extension, hence $\phi_{ij}$ is defined on both $a$ and $b$. By contradiction, the proposition follows.
\end{pf}

\begin{proposition}
  \textbf{Persistent order extension principle.} Every persistence poset $(X,\phi)$ has a linear extension.
\end{proposition}

\begin{pf}
  By proposition \ref{Zorn_cond}, conditions of Zorn's lemma are fulfilled. By Zorn's lemma $E(X,\phi)$ has maximal element $M$. Assume this element is not a linear extension of $(X,\phi)$.\\

  Then in some $M_i$, there exists incomparable pair $(a,b)$. Consider $\phi_{ij}(a)$ and $\phi_{ij}(b)$ for all $j > i$. Suppose that $\phi_{ij_1}(a) > \phi_{ij_1}(b)$ and $\phi_{ij_2}(a) < \phi_{ij_2}(b)$ for some $j_2 > j_1$. Then the map $\phi_{j_1j_2}$ cannot be defined. Hence for any $j$, the relation between images of $a$ and $b$ has the same sign if exists. If there exists such $j$ that relation between $\phi_{ij}(a)$ and $\phi_{ij}(b)$ exists we define relation between $a$ and $b$ accordingly. Otherwise, we define it arbitrarily. We can propagate the extension to the right by order preservation and to the left --- all preimages of $a$ and $b$ were incomparable and can be equipped with compatible orders.\\

  We have constructed a proper extension of $M$. By contradiction, $M$ must be a linear extension of $(X,\phi)$.
\end{pf}

\subsection{Approximation distances}

We can infer some stability results on persistence modules, using notions transferred from the category of graded modules. During this subsection, we assume conditions of Theorem 1 or Theorem 2 to be satisfied.\\

\begin{proposition} ~ \par
  \label{prop:sum} 
  Let $A$, $B$ be two persistence modules such that $d(A,0) \leq \varepsilon$ and $d(B,0) \leq \varepsilon$. Then $d(A \oplus B,0) \leq \varepsilon$.
\end{proposition}

\begin{pf}
  Recall that an index shift maps to multiplication by $m$ in a category of graded modules. We have $mA \oplus mB = m(A \oplus B)$ by definition of a direct sum. The result follows by Proposition $\ref{epstriv}$ via Theorems~1 and~2.
\end{pf}

\begin{proposition} ~ \par
  \label{prop:tensor}
  Let $A$, $B$ be two persistence modules such that $d(A,0) \leq \varepsilon$ and $d(B,0) \leq \varepsilon$. Then $d(A \otimes B,0) \leq \varepsilon$.
\end{proposition}

\begin{pf}
  The result follows from the bilinearity of the tensor product and Theorems~1 and~2.
\end{pf}

\begin{proposition}
  \label{prop:hominter}
  Let $P = \ldots \to P_n \to P_{n-1} \to \ldots$ be a persistence complex such that $P_i \stackrel{\varepsilon}{\sim} 0$ for all $i \in I$. Then the homology modules of $P$ are $\varepsilon$-interleaved with $0$.
\end{proposition}

\begin{pf}
  Assume $d_i$ is a differential in a complex. We know that $0 \to \operatorname{im}{d_{i+1}} \to \ker{d_{i}} \to H_i(P) \to 0$ is exact and that $0 \to \ker{d_i} \to P_i \to P_{i-1} \to 0$ is exact. The result follows from the application of Lemma 1 twice.
\end{pf}\\

\begin{proposition}
  Category $Fun(I, $R$-Mod)$ has enough projectives.
\end{proposition}

\begin{remark}
 This proposition is known to be true without conditions on $I$. \cite[Page 2]{Mitchell81}
\end{remark}

Since we have enough projectives, we can compute derived functors. We need the following proposition.\\

\begin{proposition}
  \label{prop:tor}
  Let $R$ be a commutative ring, $A$ and $B$ --- persistence $R$-modules such that either $A$ or $B$ is $\varepsilon$-interleaved with $0$. Then $\operatorname{Tor}_i^R(A,B) \stackrel{\varepsilon}{\sim} 0$.
\end{proposition}

\begin{pf}
  Since $R$ is commutative, $\operatorname{Tor}_i^R(A,B) = \operatorname{Tor}_i^R(B,A)$. Without loss of generality assume $B \stackrel{\varepsilon}{\sim} 0$. Let $P$ be the projective resolution of $A$. By Proposition \ref{prop:tensor}, after taking the tensor product we obtain a sequence of modules, $\varepsilon$-interleaved with $0$. The proposition follows by Proposition \ref{prop:hominter}.
\end{pf}\\

We can also derive the result about exact sequences.\\

\begin{proposition}
  \label{major}
  Let $A \xrightarrow{f} B \xrightarrow{\phi} C \xrightarrow{g} D$ be an exact sequence in the category of persistence modules. Then if $d(A,0) \leq \varepsilon$ and $d(D,0) \leq \varepsilon$, then $B \stackrel{4\varepsilon}{\sim} C$.
\end{proposition}

\begin{pf}
  In s.e.s $0 \to \ker{f} \hookrightarrow A \xrightarrow{f} \operatorname{im}f \to 0$ the $\operatorname{im}f$ is $\varepsilon$-trivial. By exactness, it is equal to $K = \ker \phi$. On the other side from $0 \to \ker{g} \xrightarrow{g} D \to D \to 0$ there follows that $d(I = \operatorname{im} \phi, 0) \leq \varepsilon$. Hence $K \stackrel{2\varepsilon}{\sim} I$ by triangle inequality.\\

  We obtain an exact sequence $0 \to K \to B \xrightarrow{\phi} C \to I \to 0$. This sequence decomposes into sequences $0 \to K \to B \to \operatorname{coIm}\phi \to 0$ and $0 \to \operatorname{Im}\phi \to C \to I \to 0$. By lemma \ref{ops} we have that $d(B,\operatorname{coIm}\phi) \leq 2\varepsilon$ and $d(C,\operatorname{Im}\phi) \leq 2\varepsilon$. Coimage and image are pointwise canonically isomorphic by the first isomorphism theorem for modules, hence $d(B, C) \leq 4\varepsilon$.
\end{pf}\\

\section{Main results}

\begin{proposition}
  \label{prop:acyclic}
  Let $A$ and $B$ be two persistence topological spaces with at least one of them being $\varepsilon$-acyclic over $R$. Then $A \star B$ is $\varepsilon$-acyclic over $R$.
\end{proposition}

\begin{pf}
  All $\operatorname{Tor}$-functors from Proposition \ref{prop:kunneth} are $\varepsilon$-interleaved with $0$ by Proposition \ref{prop:tor}. Hence by Proposition \ref{prop:sum} the right hand side of expression of Proposition \ref{prop:kunneth} is $\varepsilon$-equivalent to $0$.\\
  
  Note that on the left-hand side we operate a component-wise join. The Proposition \ref{prop:kunneth} is true component-wise and is true as a statement about persistence modules by gluing over structure maps. These structure maps are not assumed to be known, in particular, for persistent $\operatorname{Tor}_1^R(...)$.
\end{pf}\\

\begin{proposition}
  Let $x$ be an element of $(X,\phi)$. Then coverings \eqref{covering} of all components of $X$ form the persistence covering $\mathcal{U}$ with covering sets $U_1$ --- preimage of $\operatorname{st}(\mathcal{N}(x))$ under nerve functor and $U_2 = X \setminus \{x\}$.
\end{proposition}

\begin{pf}
  It suffices to check that $X \setminus \{x\}$ and preimage of $st(\mathcal{N}(x))$ are persistence subposets.
  It is evident for $X \setminus \{x\}$. Elements in the preimage of $\operatorname{st}(\mathcal{N}(x_i))$ are exactly elements comparable to $x_i$. Since structure maps preserve order, they do not move comparable elements to incomparable ones. Hence the preimage also forms a subposet.
\end{pf}\\

\begin{lemma}
  Let $(X,\phi)$ be a persistence poset and for $x=(\ldots,x_i,\ldots) \in X$ either $\mathcal{B}(X_{< x})$ or $\mathcal{B}(X_{> x})$ is $\varepsilon$-acyclic. Then persistent homology of $\mathcal{B}(X \setminus \{x\})$ and $\mathcal{B}(X)$ are $4\varepsilon$-interleaved.
\end{lemma}

\begin{pf} ~ \par
  By Proposition \ref{prop:acyclic} $\left|\operatorname{lk}(\mathcal{N}(x))\right|$ is $\varepsilon$-acyclic.\\

  Given persistence covering we can define Mayer-Vietoris exact sequence on persistence homology modules component-wise by gluing sequences for components over structure maps. Proposition \ref{major} yields the lemma.
\end{pf}

\begin{remark}
  If a map $f$ is a component-wise homotopy equivalence, it induces $0$-interleaving of homology modules.
\end{remark}

We are now ready to adapt known proof to Quillen-McCord theorem for persistence posets.

\begin{theorem} \textbf{The persistent homological Quillen-McCord theorem}\\
  Assume $X, Y$ are persistence posets of finite type indexed by very good monoid $I$, $f : X \to Y$ is an order-preserving map. Let $m$ be the number of elements of $Y$ and $R$ be a PID. Assume that the persistence topological space $\mathcal{B}(f^{-1}(Y_{\leqslant y}))$ is $\varepsilon$-acyclic over $R$ for any $y=(\ldots,y_i,\ldots) \in Y$. Then $BX$ and $BY$ are $4m\varepsilon$-interleaved over $R$.\\
\end{theorem}

\begin{pf}
  Let $\overline{X}$ be a linear extension of $X$. Recall that an element of a persistence poset has an initial element $x_i \in X_i$. We can enumerate all elements of $X$ by lexicographic order of pairs $(i,r)$ with $i$ --- index of a component in which element is born and $r$ --- the number of its initial element in order on $\overline{X_i}$.

  Let $Y^r \subset M(f)$ be the union of $Y$ and the first $r$ elements of $X$. It is a persistence subposet of $M(f)$.\\

  Consider $Y^r_{>x_r} = Y_{\geq f(x_r)}$. The persistence space $\mathcal{B}(Y_{\geq f(x_r)})$ is a component-wise cone with apex $\mathcal{B}(f(x_r))$. It is component-wise contractible, therefore $\mathcal{B}(Y^{r-1}) \hookrightarrow \mathcal{B}(Y^{r})$ is a component-wise homotopy equivalence by Lemma \ref{lem:homotopy}. By iteration, the map $\mathcal{B}(i_Y) : \mathcal{B}(Y^{0}) = \mathcal{B}(Y) \hookrightarrow \mathcal{B}(M(f)) = \mathcal{B}(Y^n)$ is a component-wise homotopy equivalence between $BY$ and $M(f)$. Note that persistence structure is not used here.\\

  Now consider a linear extension of $Y$ with the enumeration of elements and $X^r \subset M(f)$ constructed analogously. We have $X^{r-1}_{< y_r} = f^{-1}(Y_{\leqslant y_r})$. The classifying space of the latter is $\varepsilon$-acyclic over $\mathbb{F}$ by the assumption of the theorem. Hence homology modules of $\mathcal{B}(X^{r})$ and $\mathcal{B}(X^{r-1})$ are $4\varepsilon$-interleaved. By transitivity of $\varepsilon$-equivalence homology of $BX$ and $M(f)$ are $4m\varepsilon$-interleaved.\\

  We have that $H_i(BX) \stackrel{4m\varepsilon}{\sim} H_i(M(f))$ and $H_i(M(f)) \stackrel{0}{\sim} H_i(BY)$ for all $i$. Hence $H_i(BX) \stackrel{4m\varepsilon}{\sim} H_i(BY)$.
\end{pf}\\

We expect the stronger statement with no conditions on $R$ to be true with another error multiple and that it can be proved using the technique of this paper by considering error propagation in the Kunneth spectral sequence.\\

\textbf{Acknowledgements.} The author is thankful to his supervisor Anton Ayzenberg for multiple reviews of the article, invaluable support, and the link to the article \cite{Bar11} --- the starting point of the work. The author is also thankful to Ivan Limonchenko for his review.


\begin{thebibliography}{99}
  \bibitem{1}
  de Silva V, Ghrist R.
  \newblock Coordinate-free Coverage in Sensor Networks with Controlled Boundaries via Homology.
  \newblock The International Journal of Robotics Research. 2006;25(12):1205-1222.
  \newblock doi:10.1177/0278364906072252
  \bibitem{2}
  Santander, D.E. et al. (2022). 
  \newblock Nerve Theorems for Fixed Points of Neural Networks.
  \newblock In: Gasparovic, E., Robins, V., Turner, K. (eds)
  \newblock Research in Computational Topology 2. Association for Women in Mathematics Series, vol 30. Springer, Cham.
  \newblock https://doi.org/10.1007/978-3-030-95519-9\_6
  \bibitem{Kerodon002M}
  Jacob Lurie. 2023.
  \newblock Kerodon.
  \newblock https://kerodon.net/tag/002M
  \bibitem{Friedman}
  Friedman, G. (2012).
  \newblock Survey Article: An Elementary Illustrated Introduction To Simplicial Sets.
  \newblock The Rocky Mountain Journal of Mathematics, 42(2), 353–423.
  \newblock http://www.jstor.org/stable/44240054
  \bibitem{Bar11}
  Jonathan Ariel Barmak. 2011.
  \newblock On Quillen’s Theorem A for posets.
  \newblock J. Comb. Theory Ser. A 118, 8 (November 2011), 2445–2453.
  \newblock DOI:https://doi.org/10.1016/j.jcta.2011.06.008
  \bibitem{Zomorodian05}
  A. Zomorodian and G. Carlsson.
  \newblock Computing persistent homology.
  \newblock Discrete and Computational Geometry 33.2 (2005), pp. 249–274.
  \bibitem{GS16}
  Govc, Dejan \& Skraba, Primoz. (2016).
  \newblock An Approximate Nerve Theorem.
  \newblock Foundations of Computational Mathematics.
  \newblock 10.1007/s10208-017-9368-6.
  \bibitem{Corbet18}
  Corbet, R., Kerber, M.
  \newblock The representation theorem of persistence revisited and generalized.
  \newblock J Appl. and Comput. Topology 2, 1–31 (2018).
  \newblock https://doi.org/10.1007/s41468-018-0015-3
  \bibitem{Quillen72}
  Daniel Quillen,
  \newblock Higher algebraic K-theory, I: Higher K-theories Lect.
  \newblock Notes in Math. 341 (1972), 85-1
  \bibitem{Milnor56}
  Milnor, J. (1956).
  \newblock Construction of Universal Bundles, II.
  \newblock Annals of Mathematics, 63(3), second series, 430-436.
  \newblock doi:10.2307/1970012
  \bibitem{Mitchell81}
  Barry Mitchell,
  \newblock A Remark on Projectives in Functor Categories.
  \newblock Journal of Algebra 69, 24-31 (1981).
  \bibitem{Bubenik15}
  Peter Bubenik, Vin de Silva, Jonathan Scott,
  \newblock Metrics for generalized persistence modules.
  \newblock Foundations of Computational Mathematics volume 15, pages 1501–1531 (2015)
\end{thebibliography}
\end{document}